\documentclass[12pt,reqno,twoside]{amsart}
\usepackage{graphicx}
\usepackage{amsmath}
\usepackage{amssymb}

\usepackage{color}

\newcommand{\ignore}[1]{}

\newcommand{\Q}{{\mathbb {Q}}}

\newcommand{\R}{{\mathbb{R}}}  
\newcommand{\Z}{{\mathbb{Z}}}

\newcommand{\vre}{\varepsilon}

\newcommand\vare{\varepsilon}

\newcommand\AAA{\mathbb A}

\newcommand\CC{\mathbb C}

\newcommand\NN{\mathbb N}
\newcommand\RR{\mathbb R}
\newcommand\ZZ{\mathbb Z}
\newcommand\QQ{\mathbb Q}

\newcommand\cB{\mathcal{B}}

\newcommand\cD{\mathcal{D}}

\newcommand\cU{\mathcal{U}}

\newcommand\supp{\operatorname{supp}}

\newcommand\norm[1]{\left\|#1\right\|}

\newcommand\abs[1]{\left|#1\right|}

\newcommand\set[1]{\left\{{#1}\right\}}

\newtheorem{thm}{Theorem}[section]
\newtheorem{lem}[thm]{Lemma}

\newtheorem{cor}[thm]{Corollary}

\newtheorem{rem}[thm]{Remark}

\numberwithin{equation}{section}

\begin{document}

\title[Counting Diophantine approximations]{Counting intrinsic Diophantine approximations\\ in simple algebraic groups}
\author{Anish Ghosh, Alexander Gorodnik, and Amos Nevo} 
\address{Tata Institute of Fundamental research, Mumbai, India }
\email{ghosh@math.tifr.res.in}
\address{Department of Mathematics, University of Zurich }
\email{alexander.gorodnik@math.uzh.ch}
\address{Department of Mathematics, Technion IIT, Israel}
\email{anevo@tx.technion.ac.il}


\date{\today}



\begin{abstract} 
We establish an explicit asymptotic formula for the number of rational solutions of intrinsic Diophantine inequalities on simply-connected simple algebraic groups, at arbitrarily small scales. 
\end{abstract}

\maketitle

 
\section{Introduction}

\subsection{Intrinsic Diophantine approximation}

The challenge of developing the theory of Diophantine approximation
for rational points on homogeneous varieties of general algebraic groups 
was raised explicitly by S. Lang in his 1965  "Report on Diophantine Approximation" \cite[p.~189]{L65}. The key feature in the array of problems raised by Lang is that the process of Diophantine approximation called for is {\it intrinsic}, namely one seeks Diophantine approximation of the real points on the variety by the {\it rational points belonging to the variety} itself. 
Historically, investigations of distribution of rational points on homogeneous algebraic varieties focused mainly on the case of quadratic varieties, namely rational ellipsoids and hyperboloids, where several important techniques were developed and brought to bear on this problem.  Let us first mention techniques based on {\it analysis of modular forms}  e.g. in \cite{P59,M62,P79,GF85,GF87,D88,DSP90}.
Another approach based on 
{\it analysis of Hecke operators} was suggested in \cite{sar}, 
and is based on applying harmonic analysis on reductive groups
to analyze rational points on suitable homogeneous varieties.
This approach was developed, in particular, in \cite{cuo,Cl02,go,O05}.
Yet another classical approach to this problem uses the {\it circle method},
and it was undertaken in more recent times, for instance,  in \cite{HB96,M07,M14,Sar19, BVS19}.
Furthermore, the problem of Diophantine approximation on quadratic surfaces
was also studied in \cite{Dr05,KM15,FKMS, AG20} using 
{\it homogeneous dynamics} techniques, in the form of cusp excursions of diagonalizable flows on locally symmetric spaces. Finally, we refer to the survey \cite{PP17} regarding investigation of  Diophantine approximation problems in negative curvature using {\it geometric and ergodic-theoretic} techniques.

Our approach to the problem of intrinsic Diophantine approximation was initiated in \cite{GGN13,GGN14} and expanded further in \cite{GGN15,GGN18}. 
The method developed in these works 
is based on analysis of averaging operators and offers  
several advantages. In particular, it allows the consideration of rather  general homogeneous varieties, and the analysis of a wide variety of simultaneous Diophantine approximation problems in them, 
using rational points satisfying  arbitrary pre-prescribed integrality constraints.
Furthermore, the method employed is effective, making it possible to derive explicit Diophantine approximation exponents on homogeneous varieties, namely to establish a speed of approximation of a general point by a rational point on the variety. 
 To demonstrate these points, let us mention the following four natural Diophantine approximation problems arising for $\QQ$-simple linear algebraic groups $\sf G$ defined over $\QQ$:
\begin{itemize} 
\item Intrinsic Diophantine approximation on the real group variety ${\sf G}(\RR)$, using the set of {\it all rational points}. Here ${\sf G}(\RR)$ may be either compact or non-compact.

\item {\it Constrained} Diophantine approximation, involving approximation using  rational points whose denominators are divisible only by powers of primes taken from a fixed subset $S$ of primes. The set $S$ may be infinite (for example, all primes congruent to $1$ mod $4$), or finite (for example, consisting of just one prime $p$). 

\item {\it Simultaneous} approximation, where we consider the Archimedean completion together with a finite subset $S^\prime$ of non-Archimedean  completions, and approximate points on the product ${\sf G}(\RR) \times \prod_{q\in S^\prime} {\sf G}(\QQ_q)$ by all rational points, or by rational points whose denominators are constrained as above.

\item Intrinsic approximation on {\it homogeneous varieties} beyond the group varieties, where we fix an algebraic $\QQ$-subgroup ${\sf H}$ of ${\sf G}$, and consider constrained or unconstrained rational approximations on the homogeneous variety ${\sf G}/{\sf H}$, either for the variety of real points or simultaneously for several completions at once as above. 
\end{itemize}

We note that in each of these settings, two natural problems present themselves. The first is establishing  {\it uniform} Diophantine approximation with a rate of convergence valid for every single point on the variety.  
The second is the problem of {\it almost sure} Diophantine approximation, where the  rate is valid for a set of points of full measure, but is faster than the uniform rate. 
In \cite{GGN13,GGN14,GGN15} we established the existence of explicit uniform and almost sure exponents in intrinsic Diophantine approximation, as well as analogues of Khinchin's and Jarnik's theorems.  The results apply in considerable generality, namely for arbitrarily constrained intrinsic simultaneous approximation on quasi-affine homogeneous varieties of simple algebraic groups defined over an arbitrary number field. Furthermore,  in a number of cases the method provides the best possible Diophantine exponents and a sharp threshold in Khinchin's and Jarnik's theorems for the almost sure approximation problem.

\subsection{Effective density of dense subgroups}
Let us focus now on a simple linear algebraic group $\sf G$ and its group of real points ${\sf G}(\R)$, and view the previous problem from a somewhat different perspective. First note that the groups under consideration admit a rich family of countable dense subgroups.  These typically include, for example, the groups ${\sf G}(\ZZ[1/p])$ for a prime $p$, as well as
the whole group of rational points ${\sf G}(\QQ)$. Other examples of dense subgroups include ${\sf G}(\ZZ[\sqrt m])$, for suitable $m\in \NN$ which is not a perfect square, as well as  
 ${\sf G}(\QQ[\sqrt{m}])$, among many more. 
 
 It is natural to consider the problem of {\it  how dense these  subgroups actually are?}
The question of effectively quantifying the density of dense subgroups of simple compact Lie groups of the form ${\sf G}(\R)$ originated in the ground-breaking and influential analysis given by Lubotzky, Phillips and Sarnak \cite{LPS86,LPS87}, of certain dense finitely generated subgroups of the compact unitary group $\hbox{SU}_2(\CC)$. These dense countable subgroups  consist of elements in certain quaternion algebras and have matrix entries which are algebraic numbers in a fixed algebraic number field. Recently, interest in the problems stated above  has significantly increased, motivated by the goal of establishing effective bounds for the speed of topological generation of the unitary group by {\it golden gates}, namely by certain specific well-chosen generators of these (and other) dense subgroups of $\hbox{SU}_2(\CC)$. This is motivated, in turn, by the goal of estimating the efficiency of the run-time of quantum computational circuits using the given generators as the sets of gates performing the computation \cite{S15,PS18}. This problem calls specifically for {\it finitely generated} dense subgroups, and therefore highlights the significance of intrinsic Diophantine approximation with {\it constrained denominators}. For further results on $\hbox{PU}_3(\CC)$, we refer to  \cite{EP18}. 

There are several possible gauges that one can used in order to measure the density of dense subgroups. We have already mentioned the Diophantine exponent measuring the speed of approximation, or in the terminology of the previous paragraph, the speed of topological generation by the generators. This exponent has been investigated for general simple algebraic groups and their homogeneous spaces in \cite{GGN13} and more generally for lattice orbits in \cite{GGN18}. In the present paper, we will measure the density of our dense subgroups $\Gamma\subset {\sf G}(\R)$  by counting the number of elements $\gamma\in \Gamma$ of bounded size deposited in a small ball 
of radius $\delta$ centered at a general group element $x\in {\sf G}(\R)$, as $\delta \to 0^+$. We will establish an explicit asymptotic formula with an error term for this count, using as our measure of the size of elements in $\Gamma$ its algebraic "height" defined below. For the 
dense groups under consideration, this will amount to a solution of the problem of effective solution count to intrinsic Diophantine inequalities. We will consider approximation with dense groups of rational points, whose entries are {\it arbitrarily constrained}, from the finitely generated group such as ${\sf G}(\ZZ[1/p])$ for a prime $p$ to the infinitely generated group of all rational points ${\sf G}(\QQ)$.

\subsection{Main results}
\quad

{\bf Standing assumptions}. We will assume throughout the paper that ${\sf G}\subset \hbox{GL}_N$ is a simply-connected $\QQ$-almost simple linear algebraic group defined over 
$\Q$. Then the set ${\sf G}(\Q)$ of rational points is dense in the set 
${\sf G}(\R)$ of real points. More generally, if $S\subset P$ is a non-empty (finite or infinite) of the set of primes $P$, let ${\sf G}(\Z[S^{-1}])$ denotes the group consisting rational points in ${\sf G}(\RR)$ whose matrix entries all belong to the ring $\Z[S^{-1}]$ (equivalently, all entries have (reduced) denominators not divisible by primes outside $S$).  
By the Strong Approximation Property, since $\sf G$ is simply-connected, 
${\sf G}(\Z[S^{-1}])$ is dense in ${\sf G}(\R)$ provided that at least one of the groups ${\sf G}(\QQ_p)$ for $p\in S$ is non-compact. 
When this condition holds ${\sf G}$ is said to be isotropic over $ S$, and this will also be a standing assumption in our discussion. 
We refer to  \cite{PlaRa} for a detailed full discussion of strong approximation.  
\medskip

We measure the complexity of rational points
using the height function:
\begin{equation}\label{eq:height}
{\rm H}_f(r):={\prod}_{p\hbox{\tiny \,\,is prime}} \max (1,\norm{r}_p)\,,\quad r\in \hbox{Mat}_N(\Q),
\end{equation}
where $\norm{r}_p$ denote the standard maximum $p$-adic norms on $\hbox{Mat}_N(\Q_p)$.
Let us fix a left-invariant Riemannian metric $\rho$ on ${\sf G}(\R)$.
A natural measure of the quality of rational approximations 
is the following function, defined for $x\in {\sf G}(\R)$, 
$$
\omega_S(x,\delta):=\min \big\{{\rm H}_f(r):\,\, 
 \rho(x,r)\le \delta \;\;\hbox{and}\;\;
r\in {\sf G}(\Z[S^{-1}]) \big\}.
$$
The fact that  
$\omega_S(x,\delta)<\infty$ for all $x\in {\sf G}(\RR)$ and all $\delta > 0$ is equivalent to the density of 
${\sf G}(\Z[S^{-1}])$ in ${\sf G}(\R)$. Therefore, any non-trivial bound on 
$\omega_S(x,\delta)$ as $\delta\to 0^+$ quantifies the level of density of the dense group of rational approximations in question.
In \cite{GGN13} we established the following estimate on the function $\omega_S$:

\begin{thm}[\cite{GGN13}]\label{th:ggn}
	Assuming that ${\sf G}$ is isotropic over $S$, there exists an exponent $\kappa>0$ such that for all $x\in {\sf G}(\R)$,
	$$
	\omega_S(x,\delta)\le \delta^{-\kappa}\,,\quad \hbox{provided $\delta$ is sufficiently small.
	}
	$$
	The exponent $\kappa$ is given explicitly in \eqref{eq:kk} below.
\end{thm}

In the present paper, our aim is to establish a Diophantine approximation result which is far stronger. We will prove an explicit asymptotic count of the number of  constrained rational points $r$ solving the inequality $\rho(x,r)\le \delta$ at every point $x$. Furthermore, we will also establish an error estimate for this asymptotic count. 

To formulate our results explicitly, we set  for $h> \delta^{-\kappa}$,
$$
\hbox{N}_S(x,\delta,h):=\Big| \big\{r\in {\sf G}(\Z[S^{-1}]):\;  \rho(x,r)\le \delta
\;\;
\hbox{and}\;\; {\rm H}_f(r)\le h \big\} \Big|.
$$
We shall prove an explicit asymptotic formula for $\hbox{N}_S(x,\delta,h)$.

Let us write $G_p:={\sf G}(\Q_p)$ for $p\in S$, and $G_\infty:={\sf G}(\R)$.
We denote by $G_S$ the {\it restricted direct product} of $G_p$'s with $p\in S$ over the compact open subgroups $G(\widehat{\ZZ}_p)$. In the case $S=P$, $G_f:=G_P$ is the group of finite ad\'eles.
We write ${\sf G}(\mathbb{A}):=G_\infty\times G_f$ for the ad\'ele group.

 The height function $\hbox{H}_f$
extends to a proper function on $G_S$ and we set 
$$
B_S(h):=\{g\in G_S:\, \hbox{H}_f(g)\le h \}.
$$
We denote by $m_p$ the Haar measure on $G_p$ normalised so that $m_p({\sf G}(\widehat{\Z}_p))=1$
and by $m_S$ the Haar measure on $G_S$ which is the product of $m_p$'s.
The group $\Gamma_S:={\sf G}(\Z[S^{-1}])$ 
is a lattice in the product $G_\infty\times G_S$. 
We denote by $m_\infty$ the Haar measure on $G_\infty$ such that $\Gamma_S$ has covolume one in  $G_\infty\times G_S$. Given these choices, it is natural to expect that
$$
\hbox{N}_S(x,\delta,h)\approx m_\infty\big(B(x,\delta)\big)m_S \big(B_S(h)\big),
$$
where $B(x,\delta)$ denotes the Riemannian ball of radius $\delta$ centered at $x\in G_\infty$.
Indeed, this is what our main result shows:

\begin{thm}\label{the:main}
	For every $x\in G_\infty$ there exist $\theta>0$ and $c_1,c_2(x)>0$ such that
	\begin{align*}
	\hbox{\rm N}_S(x,\delta,h)= m_\infty\big(B(x,\delta)\big)m_S \big(B_S(h)\big)\Big(1+O_S\left(\delta^{-d/(d+1)} m_S(B_S(h))^{-\theta d/(d+1)} \right) \Big)
\end{align*}
where $d:=\dim({\sf G})$ and $\delta\in \big[c_1\,m_S(B_S(h))^{-\theta},c_2(x)\big)$.
Moreover, $c_2(x)$ is uniform over $x$ in bounded subsets.
\end{thm}

To give Theorem \ref{the:main} a more explicit form, let us first recall an explicit formula for the exponent $\kappa_S$ satisfying that for exponents $\kappa >  \kappa_S$ the estimate stated in Theorem \ref{th:ggn} holds.
This formula involves two parameters depending on $S$. First, for a bounded subset $\Omega \subset G_\infty$, consider the empirical rate of growth of $\Gamma_S \cap \Omega$, namely the set of $\Gamma_S$-points in $\Omega$, and define 
$$
\mathfrak{a}_S:=\sup_{\Omega\subset G_\infty} \limsup_{h\to\infty} \frac{\log |\Gamma_S\cap \Omega|}{\log h},
$$
where $\Omega$ runs over bounded subsets of $G_\infty$. For further discussion we refer to \cite[Sec. 1.2]{GGN13}.

Second, denote by $\mathfrak{q}_S$ the integrability exponent of $G_S$ in the automorphic representation, namely in the unitary representation of $G_S$ on $L^2_0({\sf G}(\mathbb{A})/{\sf G}(\Q))$, arising from the action of $G_S$ on the homogeneous space ${\sf G}(\mathbb{A})/{\sf G}(\Q)$. This exponent is discussed in detail in \cite[Sec. 3.3]{GGN13}, and we will discuss it briefly further in \S 3 below.

Then, according to \cite{GGN13} Theorem \ref{th:ggn} hold for all $\kappa>\kappa_S$, where 
\begin{equation}\label{eq:kk}
\kappa_S:=\mathfrak{q}_S d/\mathfrak{a}_S.
\end{equation}
Finally, we recall that every simple algebraic $\QQ$-group $\sf G$  determines a finite set of exceptional primes, called ramified primes, see \cite{T79} for a full account. We will exclude ramified primes from our present discussion and consider only sets $S\subset P$ for which $\sf G$ is unramified over $S$, namely $S$ contain no ramified primes for $\sf G$. 

We can now state the following special case of  Theorem \ref{the:main}, which elucidates the exponent appearing in it.  

\begin{cor}\label{cor:main}
	Assume that $\sf G$ is unramified over $S$. 
	Then for every $x\in G_\infty$, there exist $c_1,c_2(x)>0$ such that for every $\theta<(\mathfrak{q}_S d)^{-1}$,
	\begin{align*}
	\hbox{\rm N}_S(x,\delta,h)= m_\infty\big(B(x,\delta)\big)m_S \big(B_S(h)\big)\Big(1+O_{S,\theta}\left(\delta^{-d/(d+1)} m_S(B_S(h))^{-\theta d/(d+1)} \right) \Big)
	\end{align*}
	when $\delta\in \big[c_1\,m_S(B_S(h))^{-\theta},c_2(x)\big)$.
	Moreover, $c_2(x)$ is uniform over $x$ in bounded sets.
\end{cor} 

Since by \cite[Lemma~6.1]{GGN13}, 
$m_S(B_S(h))\gg_a h^{a}$ for every $a<\mathfrak{a}_S$,
Corollary \ref{cor:main} implies  that 
\begin{align*}
\hbox{\rm N}_S(x,\delta,h)= m_\infty\big(B(x,\delta)\big)m_S \big(B_S(h)\big)\Big(1+O_{S,\theta,a}\left( \delta^{-d/(d+1)} h^{-a\theta d/(d+1)}  \right) \Big).
\end{align*}
Now let us take $h=\delta^{-\kappa}$ with $\kappa>\kappa_S$.
Then $\kappa a\theta>1$ when $a<\mathfrak{a}_S$ and $\theta<(\mathfrak{q}_S d)^{-1}$ are sufficiently close to their upper bounds, and we obtain 
\begin{align*}
\hbox{\rm N}_S(x,\delta,h) &= m_\infty\big(B(x,\delta)\big)m_S \big(B_S(h)\big)\Big(1+O_{S,\theta,a}\left( \delta^{(-1 + \kappa a\theta) d/(d+1)} \right) \Big)\\
&= m_\infty\big(B(x,\delta)\big)m_S \big(B_S(h)\big) (1+o(1) )\quad\hbox{as $\delta\to 0^+$.}
\end{align*}
In particular, $\hbox{\rm N}_S(x,\delta^{-\kappa},\delta)\ne 0$
when $\delta$ is sufficiently small, and we conclude that $\omega_S(x,\delta)\le \delta^{-\kappa}$ for any $\kappa>\kappa_S$. This recovers our previous Theorem \ref{th:ggn}
in the {\it full range} of the exponents $\kappa > \kappa_S$
and moreover provides an asymptotic formula for the number of 
approximations (albeit under the unramified condition).

\begin{rem} {\rm 
\begin{enumerate}
\item We have stated Theorem \ref{the:main} and Corollary \ref{cor:main} when the family of approximating elements are taken as the intersection of $\Gamma_S$ with the growing family of height balls $B_S(h)$. However, similar results hold for many other choices of growing families,  including the height spheres, among others. This will be evident from the proof of Theorem \ref{error estimate}  in \S 2 below, on which the proofs of Theorem \ref{the:main} and Corollary \ref{cor:main} are based. General sufficient conditions for the growing family are stated in \S 2. 
%
%
%
\item 
We remark that a straightforward modification of our argument allows to establish a similar result for simultaneous approximation for products over several completions, including non-Archimedean ones, but we will not elaborate further on this case.
\end{enumerate}
}
\end{rem}
\subsection{Examples} 
We note that Theorem \ref{the:main} applies, in particular, when $\sf G$ is a simple linear algebraic group defined over $\QQ$ which is simply connected and split. In that case, the group ${\sf G}(\QQ_p)$ is non-compact for every prime $p\in P$ and so is ${\sf G}(\RR)$, namely $\sf G$ is isotropic at every place. Therefore, in this case  Theorem \ref{the:main} applies to {\it every} non-empty subset $S\subset P$, from $S=\set{p}$ to $S=P$. 

Furthermore, it is often the case that the set of ramified primes for $\sf G$ is in fact empty, and in that case the sharper form of Theorem \ref{the:main} stated in Corollary \ref{cor:main} is also valid for {\it every} non empty set $S\subset P$ without exception. 

In particular the connected simply-connected split (absolutely) almost simple linear algebraic $\QQ$-groups $\hbox{SL}_n$ satisfy both of these conditions: $\hbox{SL}_n(\QQ_p)$ is non-compact for every $p\in P$, and no prime $p$ is ramified. Corollary \ref{cor:main} quantifies the denseness of the subgroups $\Gamma_S=\hbox{SL}_n(\ZZ[S^{-1}])$ in $\hbox{SL}_n(\RR)$ by establishing the main term and an error term for the number of solutions of the intrinsic Diophantine inequalities discussed above, in the ranges specified. This applies, in particular, to the dense subgroups $\hbox{SL}_n(\ZZ[\frac1p])$ for every $p\in P$, as well as the dense subgroup  $\hbox{SL}_n(\QQ)$. 

We remark that $\hbox{SL}_2(\RR)$ is a quadratic variety, and as noted in the introduction, there exist several prior approaches to intrinsic Diophantine approximation and solution count in this case. But $\hbox{SL}_n(\RR)$ for $n \ge 3$ are not quadratic varieties, and we are not aware of any results prior to Theorem \ref{the:main} and Corollary \ref{cor:main} establishing  intrinsic Diophantine solution counts for these groups.

\subsection{Discrepancy of rational points}
Yet another perspective on the problem of counting solutions to intrinsic Diophantine inequalities arises by noting its close connection 
to the problem of estimating the discrepancy of distribution of dense sets of rational points on the varieties in question. 
 
First, recall the following general definition. Let $X$ be a locally compact metric space, and let $\nu$ be a Radon measure on $X$, positive on open sets.  Let $R_h$, $h\in \NN$ be a sequence of locally finite subsets of $X$, namely the intersection of each $R_h$ with any compact subset of $X$ is finite. Assume further that $\cup_h R_h$ is a countable dense subset of $X$. 

The collection of sets $R_h$ is said to be {\it equidistributed} in $X$ (with respect to $\nu$, and the rate function $v(h)$)  if there exists a function $v(h)\to \infty$
such that for every bounded domain $\Omega\subset X$ with boundary of measure zero
$$
\frac{|R_h\cap\Omega|}{v(h)} \to \nu(\Omega)\quad\hbox{as $h\to \infty$.}
$$
When this is the case, define the {\it discrepancy} of $R_h$ in $\Omega$ as 
\begin{equation}\label{eq:disc}
\mathcal{D}(R_h, \Omega):=\left| \frac{|R_h\cap\Omega|}{v(h)} - \nu(\Omega) \right|.
\end{equation}
Taking $\Omega=B(x,\delta)$ as metric ball in $X$, the quantity $\cD\left(R_h, B(x,\delta)\right)$ measures the pointwise discrepancy of the family of sets $R_h$ at scale $\delta$ near the point $x$. An upper bound on the function $\mathcal{D}(R_h, B(x,\delta))$ as $x$ varies in a compact subset of $X$ therefore yields a (locally uniform) pointwise upper bound on the discrepancy of the family $R_h$ at this scale. 

In the present context  we let $X={\sf G}(\RR)$, and take $R_h=R^S(h)=\Gamma_S \cap B_S(h)$. Here we view $\Gamma_S$ as a subset of ${\sf G}(\RR)$, where it is of course dense (provided ${\sf G}(\QQ_p)$ is non-compact). The set $R^S(h)$ intersects every ball 
$B(x,\delta)\subset {\sf G}(\RR)$ in a finite set, since $\Gamma_S $ is a lattice in ${\sf G}(\RR)\times G_S$, and furthermore, the union $\cup_h R^S(h)=\Gamma_S$ is dense. We take the normalization $v(h)=m_{G_S}(B_S(h))$, and the measure $\nu$ to be Haar measure on ${\sf G}(\RR)$, normalized as in the previous section. The estimates of Theorem \ref{the:main} and Corollary \ref{cor:main} are equivalent to an estimate of the expression defined in (\ref{eq:disc}). We conclude that these results establish effective equidistribution of the family $R^S(h)$, and explicit estimates of the (locally uniform)  pointwise discrepancy of the dense set of rational points $\Gamma_S $ in ${\sf G}(\RR)$.

\subsection{Method of proof}
Our arguments are  based on the general method for solving the lattice point counting problem which was developed in \cite{GN10,GN12}. Let us consider the family of compact domains 
$$
B(x,\delta)\times B_S(h)=\big\{(y,b)\in G_\infty\times G_S:\,\, \rho(y,x) \le \delta\quad\text{and}\quad {\rm H}_f(b)\le h\big\}.
$$
Then \cite[Theorem 1.9]{GN12} gives an effective solution of the lattice point counting problem for each of these domains, as $x,\delta$ is fixed and the height $h\to \infty$. The Diophantine problem under discussion amounts to producing an explicit asymptotic formula for the size of the finite sets $\abs{\Gamma_S\cap \left(B(x,\delta)\times B_S(h)\right)}$. Therefore, it can be solved provided the solution to the lattice point counting problem in question can be made effective for the family of domains varying with $\delta$. 
Hence, the required ingredient for a successful solution is bounding the error produced when the accuracy of the approximation increases, namely when the radius $\delta$ decreases to zero. 
A key feature of the method developed in \cite{GN10,GN12} is that it is based on an {\it effective mean ergodic theorem} for semisimple group actions of a form much stronger than that of the classical mean ergodic theorem. The effective mean ergodic theorem derives, in the present context, from an operator norm estimate of the averaging operators defined by the sets in question, operating in the automorphic representations. The estimate establishes decay of the operator norms as a negative power of the volume of the averaging sets, with a multiplicative constant that depends explicitly on the degree of regularity of the averaging sets. This operator norm estimate gives a bound uniform on the entire unit ball of the space of $L^2$-functions, and it is this feature  that allows the method to produce  
error estimates for the lattice point counting problem in a variable family of increasing domains, where the degree of regularity varies with $\delta$. This advantage of the method was noted in 
\cite[Remark 1.10]{GN12}, and here we develop it in detail in the present context, in order to obtain Diophantine counting results at arbitrarily small scales. 
In Section \ref{sec:lcsc} we present this argument in an abstract setting,
 and in Section  \ref{sec:finish} we apply the asymptotic formula 
 to obtained in Section  \ref{sec:lcsc} to counting rational approximations.

\subsection*{Acknowledgements}

A. Ghosh acknowledges support from a grant from the Indo-French Centre for the Promotion of Advanced Research, a Department of Science and Technology, Government of India Swarnajayanti fellowship and a MATRICS grant from the Science and Engineering Research Board. A. Ghosh acknowledges support of the Department of Atomic Energy, Government of India, under project $12-R\&D-TFR-5.01-0500$ and support from a grant from the Infosys foundation.

A. Gorodnik was supported by SNF grant 200021--182089. 
 
A. Nevo was supported by ISF Moked Grant 2919-19.

\section{Lattice points counting in variable domains}\label{sec:lcsc}

We now turn to investigate an abstract lattice point counting problem
in variable domains of general groups, which  will be the basis of our  result regarding counting solution to Diophantine inequalities. 
Let $G$ be a non-compact unimodular  locally compact second countable group, and $\Gamma\subset G$ be a discrete lattice subgroup. 
We consider a family of bounded Borel subsets $\Omega_{\delta,\cB}\subset G$ 
of positive finite measure in $G$,  
and aim to establish an asymptotic for the lattice points problem $|\Omega_{\delta, \cB}\cap \Gamma|$ uniformly. 
Our focus in the present paper will be the case when the group $G$ is a direct product of two subgroups : 
$$
G=G_\infty\times G_f \quad
\hbox{and}\quad \Omega_{\delta,\cB}:=\mathcal{O}_\delta \times \cB\,,
$$
where $G_\infty$ is an almost connected Lie group, $\cB\subset G_f$,
and where $\mathcal{O}_\delta \subset G_\infty$ are 
symmetric bounded neighborhoods of the identity in $G_\infty$. 
We will assume that $\mathcal{O}_\delta$ are decreasing as $\delta\to 0$ 
and $\cB$ is a non-empty open bounded subset of $G_f$. When applying the main counting estimate $\cB$ will be included in a growing family of sets $\cB_t$ whose volume tends to infinity. 
We write $m_{\infty}$ and $m_{f}$ for the Haar measures on the factor groups,
so that $m=m_{\infty} \times m_{f}$ is a Haar measure on $G$.
We assume that for all $0\le   r,r'\le r_0$
\begin{equation}\label{Riem-prod}
\mathcal{O}_{r} \cdot \mathcal{O}_{r'}=\mathcal{O}_{r+r'}
\quad\hbox{and}\quad 
m_\infty(\mathcal{O}_r)=r^d p(r),
\end{equation}
for fixed $d>0$ and a strictly positive Lipschitz function $p(r)$  on the interval $[0,r_0]$.  

Our main example is the case where $\mathcal{O}_r$ 
are balls of radius $r$ with respect to a left-invariant Riemannian metric on $G_\infty$, which are centered at the identity. Then \eqref{Riem-prod} follows immediately from 
invariance and standard volume computations (see, for instance, \cite[p. 
66, Cor. 5.5, Ex. 3]{Sak}).
Assumption \eqref{Riem-prod} implies the following regularity estimate for the domains $\Omega_{\delta,\cB}$.

\begin{lem}\label{l:est}
	Under the assumption \eqref{Riem-prod},
	there exists $D >  2$ such that for all $0 < \delta \le r_0/2$ and $0<\vre \le \delta/2$, 
$$
	m(\Omega_{\delta+\vre,\cB}) \le \left(1+D\frac{\vre}{\delta}\right) m(\Omega_{\delta,\cB})\quad
\hbox{and}\quad
	m(\Omega_{\delta-\vre,\cB})  \ge \left(1-D\frac{\vre}{\delta}\right) m(\Omega_{\delta,\cB}). 
$$	
The constant $D$ is given explicitly in the proof. 
\end{lem}

\begin{proof}
	Since $\delta+\vre \le r_0$, it follows from \eqref{Riem-prod} that
	$$
	\frac{m(\Omega_{\delta+\vre,\cB})}{m(\Omega_{\delta,\cB})}=\frac{m_\infty(\mathcal{O}_{\vre+\delta})}{m_\infty(\mathcal{O}_\delta)}=\frac{(\delta+\vre)^d p(\delta+\vre)}{\delta^d p(\delta)},
	$$
	and therefore 
	\begin{align*}
	\frac{m(\Omega_{\delta+\vre,,\cB})}{m(\Omega_{\delta,\cB})}-1 &\le 
	\abs{\frac{(\delta+\vre)^d \left(p(\delta+\vre)-p(\delta)\right)}{\delta^d p(\delta)}}+
	\abs{\frac{(\delta+\vre)^d}{\delta^d }-1}\\
	&\le  \left(1+\vre/\delta\right)^d \frac{|p(\delta+\vre)-p(\delta)|}{ p(\delta)}
	+(1+\vre/\delta)^d-1 \\
	&\le C_d M'M^{-1} r_0 \frac{\vre}{\delta}+2\frac{\vre}{\delta}\,,
	\end{align*}
	where $M$ is the minimum of $p$, and $M^\prime$ is the Lipschitz constant for $p$ (using also $\vre/\delta\le  1/2$).
	This proves the first estimate with a constant $D^\prime$ given by $ C_d M'M^{-1} r_0+2 > 2 $.
	Applying this estimate to the pair $(\delta-\vre,\delta)$ and using again that $\vre\le \delta/2$, we obtain
	$$m(\Omega_{\delta,\cB}) \le \left(1+D^\prime\frac{\vre}{\delta-\vare}\right) m(\Omega_{\delta-\vre,\cB})\le \left(1+2 D^\prime\frac{\vre}{\delta}\right) m(\Omega_{\delta-\vre,\cB}),$$ 
	so that 
	$$
	m(\Omega_{\delta-\vre,\cB}) \ge \left(1+2 D^\prime\frac{\vre}{\delta}\right)^{-1} m(\Omega_{\delta,\cB}) \ge \left(1-2D^\prime\frac{\vre}{\delta}\right) m(\Omega_{\delta,\cB}).
	$$
	Taking the constant $D$ to be $2D^\prime$ completes the proof.
\end{proof}
\begin{rem}
We note that the last inequality gives a positive lower bound in the range  $\vre < \delta/D$, and will be used below only under this additional assumption. 
\end{rem}
The measure $m$ on $G$ induces the invariant measure $\tilde \mu$ on the factor-space $X:=G/\Gamma$. We write $V(\Gamma):=\tilde \mu(X)$ and 
let $\mu$ denote  the unique invariant probability measure on $X$, namely 
$\mu:=V(\Gamma)^{-1}\tilde{\mu}_{X}$.
We consider the measure-preserving action of $G_f$ on the space $(X,\mu)$
and the corresponding averaging operators
$$
\pi_X(\beta)\phi(x):=\frac{1}{m_{f}(\cB)}\int_{\cB}\phi(g^{-1}x)\, dm_{f}(g),\quad \phi\in L^2(X). 
$$
We shall assume that these averages satisfy the following {\it operator norm bound:}
there exists $E(\cB)\in (0,1)$  such that
for all $\phi\in L^2(X)$,
	\begin{equation}\label{eq:norm-decay} 
	\norm{\pi_X(\beta)\phi-\int_X \phi\,d\mu}_{L^2(X)}
	\le E(\cB)\, \norm{\phi}_{L^2(X)}.
	\end{equation}

 We can now state a solution of the lattice point counting problem for a family of domains $\Omega_{\delta,\cB}$. This result constitutes a generalization of several of the main lattice point counting results in \cite{GN12} and its proof extends the methods developed there.  

\begin{thm}[counting in variable domains]\label{error estimate}
Assume $G_\infty$ is connected, $G_f$ is totally disconnected, $\mathcal{O}_\delta$ satisfy (\ref{Riem-prod}) and the operator norm bound (\ref{eq:norm-decay}) holds. Let $W\subset G_f$ be a compact open subgroup and $s,t\in G$ such that 
\begin{equation}\label{injective}
\vre_0:=\sup\{\vre>0:\; \hbox{$(\mathcal{O}_{\vre}\times W)s$ and $(\mathcal{O}_{\vre}\times W)t$ injects into $G/\Gamma$} \} > 0.
\end{equation} 
Then there exist $c_1,c_2>0$ and $A>0$ such that for every 
bounded subset $\cB$ of $G_f$ satisfying $W \cB W= \cB$,
 and every $\delta\in [c_1\,E(\cB)^{1/d},c_2)$,
\begin{align}\label{main estimate}
\abs{\frac{\abs{s\Gamma t^{-1}\cap \Omega_{\delta,\cB}}}{m(\Omega_{\delta,\cB})}-\frac{1}{V(\Gamma)}} &\le 
A\, \delta^{-d/(d+1)} E(\cB)^{-1/(d+1)}.
\end{align}
The parameters $c_1,c_2>0$ and $A>0$ are computed explicitly in the proof.
\end{thm}

\medskip
Before proceeding with the proof of Theorem \ref{error estimate}, let us note the following about its content.

\begin{rem} {\rm 
\begin{enumerate}
\item Inequality (\ref{main estimate}) does not assert that $c_1\,E(\cB)^{1/d} < c_2$, and  is not meaningful when this interval is empty. 
\item When $c_1\,E(\cB)^{1/d} < c_2$, the choice $\delta=c_1\,E(\cB)^{1/d}$ only yields an estimate which is a constant multiple of $A$, regardless of how small the spectral estimate $E(\cB)$ may be. 
\item If $\cB$ is included in a family of sets $\cB_t$ as above with $E(\cB_t)\to 0$ monotonically, then the intervals $(c_1\,E(\cB_t)^{1/d}, c_2)$
 are non-empty (for $t\ge t_0$), and approach $(0, c_2)$ as $t\to \infty$. Therefore, in this case for any fixed $0< \delta < c_2$ the limit of the estimate in the r.h.s. of (\ref{main estimate}) is zero. This establishes approximation at  arbitrarily small positive scales $\delta$. 
 \item An important special case of (3) is when the measures $m_f(\cB_t)\to \infty $, and the spectral estimate $E(\cB)$ is effective, and bounded above by a negative power of $m_f(\cB_t)$. In that case (\ref{main estimate}) establishes a simultaneous estimate, valid as $t\to \infty$ and the approximation parameter $\delta \to 0$, for a suitable range of approximation speeds. 
\end{enumerate}
}
\end{rem}

\medskip
\noindent {\it  Proof of Theorem \ref{error estimate}}. 
We consider the subsets $\mathcal{U}_\vre:=\mathcal{O}_\vre\times W$,
which provide a family of symmetric neighborhoods of the identity in $G$,
and their normalised characteristic functions
$$
\chi_\vre:=\frac{\chi_{\mathcal{U}_\vre}}{m(\mathcal{U}_\vre)}\,.
$$
Let also 
$$
\phi_\vre(h):=\sum_{\gamma\in s\Gamma t^{-1}} \chi_{\vre}(h\gamma), \quad h\in G.
$$
Since the function $\phi_\vre$ is right-invariant under 
the lattice subgroup $\Gamma^\prime:= s\Gamma s^{-1}$,
we view $\phi_\vre$ as a function on 
$G/\Gamma^\prime\simeq X$.
Clearly, $\phi_\vre$ has compact support in $X$, and 
\begin{equation}\label{eq:1}
\int_G\chi_\vre\,dm=1,\quad\quad
\int_{X}\phi_\vre\,d\tilde{\mu}=1, 
\quad\quad
\int_{X}\phi_\vre\,d\mu
=\frac{1}{V(\Gamma)}\,.
\end{equation}
Our proof is based on reduction of the 
lattice point counting problem to the operator norm estimate (and later on, the effective ergodic theorem) on $X$, combined with the regularity properties of the domain, and proceeds as follows.
We denote by $\omega_{\eta}$  the Haar-uniform averages supported on subsets $\Omega_{\eta,\cB}$ of $G$. For $\alpha > 0$ and $h\in G$,
the following estimates \eqref{eq:G_Gamma-1} and \eqref{eq:G_Gamma-2} are obviously equivalent:
\begin{align}\label{eq:G_Gamma-1}
&\abs{\pi_{X}(\omega_{\eta})\phi_\vre(h\Gamma^\prime) 
-\frac{1}{V(\Gamma)}}\le\alpha, \\
\label{eq:G_Gamma-2}
\frac{1}{V(\Gamma)}-\alpha\le &\frac{1}{m(\Omega_{\eta,\cB})}
\int_{\Omega_{\eta,\cB}}\phi_\vre(g^{-1}h\Gamma^\prime)dm(g)\le 
\frac{1}{V(\Gamma)}+\alpha\,.
\end{align}
We will estimate \eqref{eq:G_Gamma-1} using the operator norm estimate 
and Tchebysheff's inequality.   On the other hand, the integral in \eqref{eq:G_Gamma-2} is connected to counting lattice points as follows:

\begin{lem}\label{comparison}
For every $0< \vre \le \delta/4$, $\delta \le r_0/2$ and $h\in \mathcal{U}_\vre$,
\begin{equation}\label{eq:comparison}
\int_{\Omega_{\delta-2\vre,\cB}}\phi_\vre(g^{-1}h\Gamma^\prime)
\,dm(g)\le \abs{s\Gamma t^{-1}\cap \Omega_{\delta,\cB}}\le \int_{\Omega_{\delta+2\vre,\cB}}
\phi_\vre(g^{-1}h\Gamma^\prime)\,dm(g).
\end{equation}
\end{lem}

\begin{proof}
If $\chi_\vre(g^{-1}h\gamma)\ne 0$ 
for some $g\in \Omega_{\delta-2\vre,\cB}$, $h\in\mathcal{U}_\vre$, $\gamma\in s\Gamma t^{-1}$, then we obtain 
$$
\gamma\in h^{-1}\cdot \Omega_{\delta-2\vre,\cB}
\cdot \mathcal{U}_\vre\subset \Omega_{\delta,\cB}\
$$
because of \eqref{Riem-prod}. Hence, by Fubini's Theorem, the definition of $\phi_\vre$ and \eqref{eq:1},
\begin{align*}
\int_{\Omega_{\delta-2\vre,\cB}}\phi_\vre(g^{-1}h\Gamma^\prime)\,dm(g)
&= \sum_{\gamma\in s\Gamma t^{-1}} \int_{\Omega_{\delta-2\vre,\cB}}
\chi_\vre(g^{-1}h\gamma)\,dm(g)\\
&= \sum_{\gamma\in s\Gamma t^{-1}\cap \Omega_{\delta,\cB}} \int_{\Omega_{\delta-2\vre,\cB}}
 \chi_\vre(g^{-1}h\gamma)\,dm(g)\le \abs{s\Gamma t^{-1}\cap \Omega_{\delta,\cB}}.
\end{align*}
In the other direction, 
for $\gamma\in s\Gamma t^{-1}\cap \Omega_{\delta,\cB}$ and $h\in\mathcal{U}_\vre$, using \eqref{Riem-prod}
$$
\supp(g\mapsto \chi_\vre(g^{-1}h\gamma))
=h\gamma(\supp\,\chi_\vre)^{-1}\subset \Omega_{\delta+2\vre,\cB}.
$$
Since $\chi_\vre\ge 0$ and \eqref{eq:1} holds, we conclude that 
\begin{align*}
\int_{\Omega_{\delta+2\vre,\cB}}\phi_\vre(g^{-1}h\Gamma^\prime)\,dm(g)
&= \sum_{\gamma\in s\Gamma t^{-1}} \int_{\Omega_{\delta+2\vre,\cB}}\chi_\vre(g^{-1}h\gamma)\,dm(g)\\
&\ge \sum_{\gamma\in s\Gamma t^{-1}\cap \Omega_{\delta,\cB}}
 \int_{\Omega_{\delta+2\vre,\cB}}\chi_\vre(g^{-1}h\gamma)\,dm(g)= \abs{s\Gamma t^{-1}\cap \Omega_{\delta,\cB}},
\end{align*}
as required.
\end{proof}

The Haar-uniform measures $\omega_{\eta}$ on the sets $\Omega_{\eta,\cB}$ are the product of the Haar-uniform measures $\omega_\eta^\infty$ supported on $\mathcal{O}_\eta$ and the Haar-uniform measure $\beta$ on $\cB$,
so that 
$$
\pi_{X}(\omega_{\eta})=\pi_X(\omega^\infty_\eta)\pi_X(\beta).
$$
Moreover, the operators $\pi_X(\omega^\infty_\eta)$ and $\pi_X(\beta)$ commute.
Therefore, we deduce from \eqref{eq:norm-decay} that
\begin{align*}
\norm{\pi_{X}(\omega_{\eta})\phi_\vre-\int_{X} \phi_\vre\, d\mu }_{L^2(X)} &\le \norm{\pi_{X}(\omega_\eta^\infty)}\norm{\pi_{X}(\beta) \phi_\vre-\int_{X} \phi_\vre\, d\mu}_{L^2(X)} \\
&\le E(\cB)\norm{\phi_\vre}_{L^2(X)}\,,
\end{align*}
and by Tchebysheff's inequality for all $\alpha>0$,
\begin{align*}
\mu\left(\set{h\Gamma^\prime\,:\, \abs{\pi_{X}
(\omega_{\eta})\phi_\vre(h\Gamma^\prime) -\frac{1}{V(\Gamma)}} > \alpha}\right)
\le 
\alpha^{-2}  E(\cB)^{2}
\norm{\phi_\vre}^2_{L^2(X)}.
\end{align*}

We shall additionally assume the parameter $\vre$ satisfies  $\vre<\vre_0/2$.
Then the projection $\mathcal{U}_{2\vre}\to \mathcal{U}_{2\vre}s\Gamma$
is injective. Since the neighborhoods $\mathcal{U}_\vre$ are symmetric and 
$\mathcal{U}_{2\vre}=\mathcal{U}_{\vre}^2$, we also obtain that
\begin{equation}\label{eq:m0}
\mathcal{U}_{\vre}\rho_1\cap \mathcal{U}_{\vre}\rho_2=\emptyset\quad
\hbox{for all $\rho_1\ne \rho_2\in\Gamma^\prime$,}
\end{equation}
in particular,
\begin{equation}\label{eq:m1}
\mu(\mathcal{U}_\vre\Gamma^\prime)=\frac{m(\mathcal{U}_\vre)}{V(\Gamma)}.
\end{equation}
We observe that
\begin{align*}
\norm{\phi_\vre}^2_{L^2(X)}&=\int_{X}
\left({\sum}_{\gamma\in s\Gamma t^{-1}} \chi_\vre(h\gamma)\right)^2  d\mu(h\Gamma^\prime)\\
&=
\int_{X}
\left({\sum}_{\gamma_1,\gamma_2\in s\Gamma t^{-1}} \chi_\vre(h\gamma_1)\chi_\vre(h\gamma_2) \right)  d\mu(h\Gamma^\prime).
\end{align*}
If $h\gamma_1, h\gamma_2\in \supp(\chi_\vre)=\cU_\vre$,
then $\gamma_2^{-1}\gamma_1\in t\Gamma t^{-1}\cap \cU_\vre^2$,
but since $\vre<\vre_0/2$, the map $\cU_\vre^2=\cU_{2\vre}\to \cU_{2\vre}t\Gamma$ is injective. Therefore, we conclude that $\gamma_1=\gamma_2$, and 
\begin{align*}
\norm{\phi_\vre}^2_{L^2(X)}
= \int_X \left({\sum}_{\gamma\in s\Gamma t^{-1}} \frac{\chi_{\cU_\vre}(h\gamma)}{m(\cU_\vre)^2}\right)\, d\mu(h\Gamma^\prime)=
\frac{m(\mathcal{U}_\vre)^{-1}}{V(\Gamma)},
\end{align*}
so that
\begin{equation}\label{eq:measure}
\mu\left(\set{h\Gamma\,:\, \left|\pi_{X}
(\omega_{\eta})\phi_\vre(h\Gamma^\prime) -\frac{1}{V(\Gamma)}\right| > \alpha}\right)
\le 
  \alpha^{-2}   E(\cB)^{2} \frac{m(\mathcal{U}_\vre)^{-1}}{ V(\Gamma)} \,.
\end{equation}
When the parameter $\alpha$ is sufficiently large, we obtain
\begin{equation}\label{eq:intersection}
\mathcal{U}_\vre\Gamma^\prime\cap 
 \set{h\Gamma^\prime\,:\, \abs{\pi_{X}
(\omega_{\eta})\phi_\vre(h\Gamma^\prime) -\frac{1}{V(\Gamma)}} \le \alpha}
\neq \emptyset.
\end{equation}
Let us now write $\eta=\delta+2\vre$ with $\delta\le r_0/2$ and $\vre<\delta/4$, so that averages $\omega_{\eta}$ under consideration  are supported on the sets $\Omega_{\delta+2\vre,\cB}$. Applying (\ref{eq:G_Gamma-1})
 and (\ref{eq:G_Gamma-2}) to any $h$ in the non-empty intersection 
 (\ref{eq:intersection}), we obtain that 
$$
 \frac{1}{m(\Omega_{\delta+2\vre,\cB})}
\int_{\Omega_{\delta+2\vre,\cB}}\phi_\vre(g^{-1}h\Gamma^\prime)dm(g)\le 
\frac{1}{V(\Gamma)}+\alpha.
$$
On the other hand, by Lemma \ref{comparison},  for any $h\in \mathcal{U}_\vre$, 
$$
\abs{s\Gamma t^{-1}\cap \Omega_{\delta,\cB}}\le 
\int_{\Omega_{\delta+2\vre,\cB}}\phi_\vre(g^{-1}h\Gamma^\prime)dm(g)\,\,.
$$
Combining these estimates, we conclude that 
$$
\abs{s\Gamma t^{-1}\cap \Omega_{\delta,\cB}}\le
 \left(\frac{1}{V(\Gamma)}+\alpha\right)m(\Omega_{\delta+2\vre,\cB})\,.
$$
Since $0 < \delta \le r_0/2$ and $0 < \vre\le \delta/4$,
it follows from Lemma \ref{l:est} that
$$
m(\Omega_{\delta+2\vre,\cB})\le \left(1+2D\frac{\vre}{\delta}\right)m(\Omega_{\delta,\cB}),
$$
and therefore 
\begin{equation}\label{eq:bbbb0} 
\abs{s\Gamma t^{-1}\cap \Omega_{\delta,\cB}}\le\left(\frac{1}{V(\Gamma)}+\alpha\right) \left(1+2D\frac{\vre}{\delta}\right)m(\Omega_{\delta,\cB}).
\end{equation} 
This inequality holds as soon as  (\ref{eq:intersection}) holds, 
and so certainly if we have 
\begin{equation}\label{eq:neq}
  \alpha^{-2}   E(\cB)^{2} \frac{m(\mathcal{U}_\vre)^{-1}}{ V(\Gamma)}\le 
\frac{1}{4} \cdot\frac{m(\mathcal{U}_\vre)}{V(\Gamma)}.
\end{equation}
Indeed, then the right hand side of (\ref{eq:measure}) is strictly 
smaller than $\mu
(\mathcal{U}_\vre^\prime)=m(\mathcal{U}_\vre)/V(\Gamma)$, so that the intersection 
(\ref{eq:intersection}) is necessarily non-empty. 
We set 
$$\alpha=2 m(\mathcal{U}_\vre)^{-1} E(\cB),$$
 so that the inequality (\ref{eq:neq}) holds. 
 
 We now additionally assume that $\vre <  \delta/(2D)$.
Then we deduce from \eqref{eq:bbbb0} that 
 \begin{align*}
\frac{\abs{s\Gamma t^{-1}\cap \Omega_{\delta,\cB}}}{m(\Omega_{\delta,\cB})}-\frac{1}{V(\Gamma)} 
\le \alpha\left(1+ \frac{2D\vre}{\delta}\right)+\frac{2D}{V(\Gamma)}\frac{\vre}{\delta}\le 2\alpha+\frac{2D}{V(\Gamma)}\frac{\vre}{\delta}\,.
\end{align*}

As for the lower bound in Theorem \ref{error estimate}, we note that applying (\ref{eq:G_Gamma-1})
 and (\ref{eq:G_Gamma-2}) to any $h$ in the non-empty intersection 
 (\ref{eq:intersection}), with the choice $\eta=\delta-2\vre$
 with $\delta\le r_0/2$ and $\vre\le \delta/4$, we obtain
$$
 \frac{1}{m(\Omega_{\delta-2\vre,\cB})}
\int_{\Omega_{\delta-2\vre,\cB}}\phi_\vre(g^{-1}h\Gamma^\prime)dm(g)\ge 
\frac{1}{V(\Gamma)}-\alpha,
$$
and furthermore,  by Lemma \ref{comparison}, 
$$
\abs{s\Gamma t^{-1}\cap \Omega_{\delta,\cB}}\ge 
\int_{\Omega_{\delta-2\vre,\cB}}\phi_\vre(g^{-1}h\Gamma^\prime)dm(g)\,\,.
$$
Combining these estimates, we conclude that 
$$
\frac{\abs{s\Gamma t^{-1}\cap \Omega_{\delta,\cB}}}{ m(\Omega_{\delta-2\vre,\cB})}\ge
 \frac{1}{V(\Gamma)}-\alpha   \,.
$$
According to Lemma \ref{l:est}, when $\frac{D\vre}{\delta}<  1$
$$
m(\Omega_{\delta-\vre,\cB})\ge \left(1-D\frac{\vre}{\delta}\right)m(\Omega_{\delta,\cB}) > 0\,,
$$
so that when $\frac{2D\vre}{\delta}<  1$, we can conclude 
$$
\frac{\abs{s\Gamma t^{-1}\cap \Omega_{\delta,\cB}}}{\left(1-2D\frac{\vre}{\delta}\right)m(\Omega_{\delta,\cB})}\ge\frac{1}{V(\Gamma)}-\alpha\,.
$$
Therefore, 
$$
\frac{\abs{s\Gamma t^{-1}\cap \Omega_{\delta,\cB}}}{m(\Omega_{\delta,\cB})}\ge\left(\frac{1}{V(\Gamma)}-\alpha\right) \left(1-2D\frac{\vre}{\delta}\right)\,.
$$
and 
$$
\frac{\abs{s\Gamma t^{-1}\cap \Omega_{\delta,\cB}}}{m(\Omega_{\delta,\cB})}-\frac{1}{V(\Gamma)} 
\ge  -\alpha\left(1-\frac{2D\vre}{\delta}\right)-\frac{2D}{V(\Gamma)}\frac{\vre}{\delta}\ge
-2\alpha-\frac{2D}{V(\Gamma)}\frac{\vre}{\delta}\,. 
$$
Combining the established upper and lower bounds, we obtain
\begin{align*}
\abs{\frac{\abs{s\Gamma t^{-1}\cap \Omega_{\delta,\cB}}}{m(\Omega_{\delta,\cB})}-\frac{1}{V(\Gamma)}} &\le 
2\alpha+\frac{2D}{V(\Gamma)}\frac{\vre}{\delta}\\
&=4m(\mathcal{U}_\vre)^{-1} E(\cB)+\frac{2D}{V(\Gamma)}\frac{\vre}{\delta}\\
&=4M^{-1} m_f(W)^{-1} \vre^{-d} E(\cB)+\frac{2D}{V(\Gamma)}\frac{\vre}{\delta},
\end{align*}
where $M$ is the minimum of the function $p$ from \eqref{Riem-prod}.
Choosing the parameter $\vre$ as
$$
\vre=\left(M^{-1} D^{-1}\, m_f(W)^{-1}\, E(\cB)\, V(\Gamma)  \delta \right)^{1/(d+1)},
$$
we obtain the estimate
\begin{align*}
\abs{\frac{\abs{s \Gamma t^{-1}\cap \Omega_{\delta,\cB}}}{m(\Omega_{\delta,\cB})}-\frac{1}{V(\Gamma)}} &\le 
A\, \delta^{-d/(d+1)} E(\cB)^{1/(d+1)}
\end{align*}
with
\begin{equation}\label{eq:A}
A:=6 M^{-1/(d+1)}D^{d/(d+1)}\, m_f(W)^{-1/(d+1)}  V(\Gamma)^{-d/(d+1)}.
\end{equation}
It remain to insure that our choice of the parameter $\vre$ satisfies 
the conditions which we have imposed. We required that
$$
\vre<\vre_0/2\quad\hbox{and}\quad \vre < \delta/(2D).
$$
These conditions are satisfied provided that $c_1 E(\cB)^{1/d}\le \delta <c^\prime_2$, where
\begin{equation}\label{eq:c1}
c_1:= 2^{(d+1)/d}D M^{-1/d} \, m_f(W)^{-1/d}V(\Gamma)^{1/d},
\end{equation}
and 
\begin{equation}\label{eq:c2}
c^\prime_2:= 2^{-(d+1)} M D\,m_f(W) V(\Gamma) \vre_0^{d+1}.
\end{equation}
Here we used that $E(\cB)\in (0,1]$. Finally, we set $c_2:=\min\{c^\prime_2,r_0/2\}$.
 \qed

\section{Proof of the main results}\label{sec:finish}

Before turning to the proof of the main results, let us note the following 
about the automorphic representation, namely the unitary representation 
of the ad\'ele group ${\sf  G}(\AAA)$ on the space $L^2({\sf  G}(\AAA)/{\sf G}(\QQ))$, arising from the action of ${\sf  G}(\AAA)$ by translation 
on the space ${\sf  G}(\AAA)/{\sf G}(\QQ)$ equipped with invariant probability measure. Denote the automorphic representation restricted to the subgroup $G_p={\sf G}(\QQ_p)$ by $\pi_p$, and denote by $\pi_p^{0}$ its further restriction to the subspace of zero-integral functions $L_0^2({\sf  G}(\AAA)/{\sf G}(\QQ))$.

When $\sf G$ is simply connected, the {\it integrability exponent} of $\pi_p$ w.r.t. the compact open subgroup ${\sf G}(\widehat{\ZZ}_p)\subset \QQ_p$ is defined,  for primes $p$ which are unramified and for which $\sf G$ is isotropic, as follows:
\begin{equation}\label{eq:q_v}
\mathfrak{q}_p({\sf G}):=\inf\left\{ q \ge 2:\, \begin{tabular}{l}
$\forall$\hbox{ ${\sf G}(\widehat{\ZZ}_p)$-inv. $w\in L_0^2({\sf  G}(\AAA)/{\sf G}(\QQ))$}\\
\hbox{$\left< \pi_p(g)w,w\right>\in L^q({\sf G}(\QQ_p)$})
\end{tabular}
\right\}.
\end{equation}
 An equivalent definition is that 
$\mathfrak{q}_p({\sf G})$ is the infimum of  $q\ge 2$ having the property that every irreducible ${\sf G}(\widehat{\ZZ}_p)$-spherical unitary representation  
$ \tau_p$ of ${\sf G}(\QQ_p)$ which is weakly contained in $ \pi_p$ is an  $L^{\mathfrak{q}+}$-representation, namely has a dense subspace 
giving rise to matrix coefficients in $L^q({\sf G}(\QQ_p))$, for every $q >\mathfrak{q}_p({\sf G})$. We refer to \cite[\S 3.3]{GGN13} for further discussion, and also for definition of the integrability exponent for a general subset $S\subset P$. 
 
 We note that when $\sf G$ is not simply connected the definition above must be modified, by changing the subspace $L_0^2({\sf  G}(\AAA)/{\sf G}(\QQ))$ to the subspace orthogonal to all the automorphic characters. When the prime $p$ is ramified, again the definition above must be modified, by changing from ${\sf G}(\widehat{\ZZ}_p)$ to a possibly different compact open subgroup.  We will not elaborate on these matters further and instead restrict our discussion accordingly, namely to simply-connected groups and unramified primes.

It is well-known that the unitary representation of $G_\infty\times G_S$  on 
$L^2_0\left((G_\infty\times G_S)/\Gamma_S\right)$ is equivalent to the unitary representation of $G_\infty\times G_S$ on a suitable subspace of $  L_0^2({\sf  G}(\AAA)/{\sf G}(\QQ))$. As a result, 
an upper bound for the integrability exponent as defined is valid also for the unitary representation which underlies the discussion in \S 2, namely the unitary representation of $G_\infty\times G_S$  on $L^2_0\left((G_\infty\times G_S)/\Gamma_S\right)$.

\begin{proof}[Proof of Theorem \ref{the:main} and Corollary \ref{cor:main}]
We use Theorem \ref{error estimate} for the group $G_\infty\times G_S$
and the family of sets 
$\Omega_{\delta,\cB}:=B(e,\delta)\times B_S(h)$.
The condition \eqref{Riem-prod} follows from the properties of Riemannian balls stated in \cite[p. 
66, Cor. 5.5, Ex. 3]{Sak}, and the sets $B_S(h)$ are bi-invariant under the compact open subgroup 
$W:=\prod_{p\in S} {\sf G}(\widehat \Z_p)$. It follows from the invariance
of the metric that $B(x,\delta)=xB(e,\delta)$. Therefore,
$$
\hbox{N}_S(x,\delta,h)=\abs{\Gamma_S\cap B(x,\delta)\times B_S(h)}=
\abs{x^{-1}\Gamma_S\cap B(e,\delta)\times B_S(h)}.
$$
It remains to verify the operator norm condition \eqref{eq:norm-decay}
for the Haar-uniform averages $\pi_X(\beta_h)$ supported on the sets $B_S(h)$, acting on the space $X=(G_\infty\times G_S)/\Gamma_S$. An effective mean ergodic theorem 
for actions of simply-connected $\QQ$-simple ad\'ele groups ${\sf G}(\mathbb{A})$, valid for these operators, was established in \cite[Cor. 6.7]{GN12}. The same argument applies to the group $G_S$, and holds for any action of $G_S$ satisfying the integrability condition. 

When $S$ is finite, the integrability condition is the content of property
$(\tau)$ of automorphic representations established in \cite{Cl03}.
For infinite $S$, the integrability also holds as was verified by \cite[Th.~3.20 and Th.~3.7]{GMO}. We obtain that  there exist $c,\tau> 0$
such that for every $\phi\in L^2(X)$,
\begin{equation}\label{eq:ergg}
\norm{\pi_X(\beta_h)\phi-\int_X \phi\,d\mu}_{L^2(X)}
\le \min\big(1, c\, m_S(B_S(h))^{-\tau}\big)\, \norm{\phi}_{L^2(X)}.
\end{equation}
Hence, Theorem \ref{error estimate} can be applied in our setting.
Moreover, we note that in the constants $A,c_1,c_2$ given by \eqref{eq:A}, \eqref{eq:c1}, \eqref{eq:c2}, the only dependence on $x$
comes from the injectivity radius $\vre_0$, which appears only in $c_2$,
and $\vre_0$ is uniform over bounded sets.
Hence, we conclude that Theorem \ref{the:main} holds with $\theta=\tau/d$.

Regarding Corollary \ref{cor:main}, we note that  
\cite[Th. 4.2]{GGN13} implies that in the unramified case the ergodic theorem holds with $\tau <\mathfrak{q}_S^{-1}$, so that in this case the estimate of Theorem \ref{the:main} holds when $\theta< (\mathfrak{q}_S d)^{-1}$.
\end{proof}

\end{document}